\documentclass[12pt]{article}
\usepackage{graphicx}
\usepackage{multicol}
\usepackage{color}

\usepackage{geometry}
\begin{document}

\begin{center}

{ \bf 
Unbounded eigenfunctions  in the stability problem

for a three-layer flow in porous media}

 \vspace{0.25cm}
 
 Gelu Pa\c{s}a

Simion Stoilow Institute
of Mathematics, 
Calea Grivitei 21, 
Bucharest, Romania

{\it  E-mail}:  gelu.pasa@imar.ro  
\hspace{0.5cm}   and
\hspace{0.5cm}
pasa.gelu@gmail.com

\end{center}

\vspace{0.25cm}

\noindent
We study the linear stability of the  displacement of three
Stokes fluids with 
constant viscosity 
in a porous medium when
the middle   fluid is contained in a bounded
region.
We use the Hele-Shaw model.
The 
eigenfunctions of the stability system are the
amplitudes
of the linear perturbations.  
These amplitudes  must be small.
We get  unbounded eigenfunctions.
So the stability problem 
has   no physical sense.

 \vspace{0.25cm}




\begin{center}
{\bf 
1. Introduction}
\end{center}

We consider   the linear stability 
problem for the displacement of 
three Stokes fluids
in a Hele-Shaw cell, parallel with the horizontal fix plane $x_1Oy$ - see \cite{BE}, 
\cite{SAFF-TAY}. The three constant viscosities 
$\mu_L, \mu, \mu_R$
are (in fact) divided by the constant permeability of the equivalent porous medium.
The fluid
$\mu_L$ is displacing the middle fluid 
$\mu$, which in turn is pushing a fluid with viscosity $\mu_R$, 
in the positive direction
$Ox_1$.
The  middle liquid  $\mu$ is contained in a bounded region. We suppose
\begin{equation}\label{U}
\mu_L < \mu < \mu_R.
\end{equation}

A 
large number of successive
intermediate liquids can be considered between the initial displacing fluids
with viscosities $\mu_L, \mu_R$.
On this way some authors
try to minimize the   Saffman-Taylor instability
(see \cite{DA1} - \cite{DA6}), by using the seminal paper  \cite{DA0}.
The obtained 
results have some weak points, explained in some works of the present author (see
\cite{PA00} -\cite{PA4}).
We give here 
only a selective list of
related papers. Some have been published recently and contain more 
detailed references.

The linear stability system was obtained in \cite{DA0} 
and also in \cite{PA3}.
An  existence condition  for the non-zero
solution of the stability system was used in \cite{DA0} to obtain a
quadratic equation for the  the growth rates; some numerical results were given.
A compatibility condition
was derived in \cite{PA3}; we proved that this condition is not fulfilled. Thus the stability problem has no solution. However, some (natural) additional hypothesis were used in 
\cite{PA3}. 

 In this paper we use the above existence
condition  for  obtain the behavior of the
 eigenfunctions 
of the stability system 
in the range of very large wavenumbers. 
The new element of this paper is following.
 We get
 unbounded eigenfunctions for  increasing  wavenumbers.
Therefore the perturbations are also
unbounded for large wavenumbers.
Thus the stability problem 
has   no physical sense. This time, we not use any additional hypothesis.

\begin{center}
{\bf 2. The stability 
problem}
\end{center}

The displacement  is given
 by the velocity $(U,0)$ of the 
 fluid $\mu_L$ far upstream.
We use  the moving reference $x=x_1-Ut$,
where $t$ is the time.
 The intermediate liquid is contained in the segment $(a,b), b 
\leq 0$.
The relevant (basic) flow equations are  quite similar with
 Darcy's law in porous media.
A stationary solution
exists,  with
two straight  interfaces 
$x=a, \,\, x=b$
between the  fluid layers -
see \cite{GOR-HOM-1}.
The Laplace-Young law
acts 
on both  interfaces, where   two positive  surface tensions 
$T(a),T(b)$ exist. 
 
We perform a linearized stability analysis by normal modes. 
We use the Fourier decomposition of the velocity perturbation
 in the $Ox$ direction: 
\begin{equation}\label{J1}
 u'= {\color{red} \epsilon} 
f(x)\exp(iky + 
\sigma t).
\end{equation}                      
Here 
  $k, \sigma, f$ are 
 the wave numbers, growth rates and 
  amplitudes. 
  $\epsilon $ is a small positive number.
 The function $f$ must
  be bounded. Moreover,  $f$ must 
 decay to zero in the far field and is continuous.
 The derivative of $f$
 is not continuous  in $a,b$. In 
 \cite{SAFF-TAY}
are  imposed the conditions 
$$
f(x)= f(a)e^{k(x-a)}, \,\, \forall x \leq a,
\,\, k\geq 0 ; $$
\begin{equation}\label{J2}
f(x)= f(b)e^{-k(x-b)}, \,\, \forall x \geq b,
\,\, k \geq 0 .
\end{equation}
We insert (\ref{J1}) in
the linearized disturbance equations obtained from the Darcy's law.
We follow the procedure
given  in  \cite{DP}, \cite{DA0} and  we 
obtain
the problem (2a), (2b), (3)  of  \cite{DA0}:
 \begin{equation}\label{A5}
 f_{xx}(x,k)-k^2 f(x,k)=0, \,\,
 x  \in (a,b), \quad 
 b \leq 0, \quad
 \forall k \geq 0;
\end{equation}
\begin{equation}\label{A6}
 f_x^+(a,k)= 
 [-k E_a(k)/ \sigma(k) + 
 \mu_L k/ \mu]f(a,k);
\end{equation}
\begin{equation}\label{A7}
 f_x^-(b,k)=
[k E_b(k) / \sigma(k) - \mu_R k/ \mu ]f(b,k);
\end{equation}
$$
 E_a(k): =  \frac{
 (\mu-\mu_L)Uk
 -  T(a)k^3}{\mu} ,
                $$
 \begin{equation}\label{A9}
 E_b(k):=\frac{
 (\mu_R-\mu)Uk-
 T(b)k^3}{\mu} ;
\end{equation}
$f_x$ is the derivative; 
$f_x^-, f_x^+$ are 
 the lateral 
 derivatives; $f$ is depending also on $k$.

\begin{center}
{\bf 3. The non-existence 
result}
\end{center}

The general solution of 
(\ref{A5}) is
\begin{equation}\label{B2}
f(x,k)=A(k)e^{kx}+
B(k)e^{-kx}  .         
\end{equation}
The coefficients $A,B$ were considered constant in \cite{DA0}. 
We use the relations 
(\ref{A6})-(\ref{A7}) and obtain a stability system 
for $A(k), B(k)$, 
with variable coefficients 
$c(k),d(k),g(k),h(k)$  (see formulas 
(2b), (3), (5), pag. 112101-5 in \cite{DA0}
and the system (12) of
\cite{PA3}):
\begin{equation}\label{B3}
 A(k)e^{ka} c(k) + 
 B(k)e^{-ka} d(k)=0,
 \quad        
 A(k)e^{kb} g(k) + 
 B(k)e^{-kb} h(k)=0;
\end{equation}    
 \begin{equation}\label{B4}  
c(k)=[\mu_L-\mu-E_a(k)/\sigma(k)],\quad 
d(k)=[\mu_L+\mu-E_a(k)/\sigma(k)],
\end{equation}
\begin{equation}\label{B5}
g(k)=[\mu_R+\mu-E_b(k)/\sigma(k)],\quad 
h(k)=[\mu_R-\mu-E_b(k)/\sigma(k)].
\end{equation}

A nontrivial solution 
$(A(k),B(k))$ exists if the following condition is verified 
\begin{equation}\label{B5A}
 e^{2k(a-b)}c(k)h(k)- g(k)d(k)=0,
 \quad \forall k \geq 0.
\end{equation}

 Our main goal is to find
$\lim_{k\to\infty}A(k)$ and
$\lim_{k\to\infty}B(k)$
by using (\ref{B3}) - 
(\ref{B5A}).

\vspace{0.5cm}

{\it Remark 1}. We have 
\begin{equation}\label{P1}
d(k),c(k),g(k),h(k)\neq 0, \forall k   > 0.
\end{equation}
Indeed, suppose
$\exists \,\,
 m \,\, s.t. \,\, 
d(m)=0, c(m)= -2\mu$. 
From $(\ref{B3})_1$  and 
(\ref{B2})
we get  
$$ A(m)=0, \,
f(x, m)=B(m)e^{-mx},  \, 
f_x(x,m)=-mB(m)e^{-mx}.
                      $$ 
 We emphasize that 
$B(m)\neq 0$, otherwise
$f(x,m) \equiv 0$; the eigenfunctions can not be 
identically zero.  
We simplify with 
$B(m)$,  then   
(\ref{A6})-(\ref{A7}) in
$k=m$ lead us to
\begin{equation}\label{B9}
\mu_L e^{-ma}+
\mu e^{-ma} =      
\frac{E_a(m)}{\sigma(m)}e^{-ma};       \quad 
\mu_R e^{-mb}+
\mu e^{-mb} = 
\frac{E_b(m)}{\sigma(m)}
e^{-mb}.
\end{equation}
In both relation we have the same $\sigma(m)$ then  we get 
$ E_a(m)/(\mu_L+\mu)=
E_b(m)/(\mu_R+\mu).     $
This is an unexpected restriction on $\mu$,
 not considered  
as a hypothesis.

We get
 $ c(k), g(k), h(k) \neq 0, \forall \, k >0 $ by using the same arguments.
For example, suppose
$\exists \,\,
 n \,\, s.t. \,\, 
c(n)=0$. Then  (\ref{B2}), 
$(\ref{B3})_1$  and 
(\ref{A6})-(\ref{A7}) in
$k=n$   give us 
$ B(n)=0, \,\, 
f(x,n)= A(n)e^{nx}, \,\,
A(n)\neq0, \,\,
f_x(x,n)=nA(n)e^{nx} $,
thus $E_a(n)/(\mu_L-\mu)=
E_b(n)/(\mu_R+\mu)$.

\hfill $\bullet $

\vspace{0.25cm}

{\it Remark 2.} 
The following properties hold:
\begin{equation}\label{P2}
\lim_{k\to\infty} d(k)
=0 \Rightarrow 
\lim_{k\to\infty} g(k),
 \lim_{k\to\infty} c(k),
 \lim_{k\to\infty} h(k)
\quad \mbox {are real bounded (not zero) quantities}.
 \end{equation}
\begin{equation}\label{P3} 
\lim_{k\to\infty} g(k)
=0 \Rightarrow 
\lim_{k\to\infty} d(k),
 \lim_{k\to\infty} c(k),
 \lim_{k\to\infty} h(k)
 \quad \mbox {are real bounded (not zero) quantities}.
 \end{equation}
 Indeed, suppose
 $\lim_{k\to\infty}d(k)=0$ 
 . Thus 
 $ \lim_{k\to\infty} 
 E_a/\sigma = 
 (\mu_L+\mu)        $
 and we have 
 $$\lim_{k\to\infty} g(k)=
 \mu_R+\mu- (T_b/T_a)
 (\mu_L+\mu),         $$
 $$\lim_{k\to\infty} h(k)=
 \mu_R-\mu- (T_b/T_a)
 (\mu_L+\mu),      \quad
 \lim_{k\to\infty} c(k)=
 - 2 \mu.  $$
 We use the  same arguments and  get  (\ref{P3}).

 \hfill   $\bullet$ 
  
 \vspace{0.25cm}

 For a function $H(k)$,
 we use also the notation
 $H(k)\rightarrow \alpha$
 when $\lim_{k\to\infty}H(k)=
 \alpha$. 
 
 \hfill $\bullet$
 
 \vspace{0.5cm}

 From {\it Remark 2} it follows that $d(k)\rightarrow 0 $ is a 
 possible ``solution'' of 
 (\ref{B5A}) when 
 $k \rightarrow \infty$.
 Indeed, the limits of $c(k),h(k),g(k)$ are bounded, thus
 $e^{2k(a-b)}c(k)h(k)
 \rightarrow 0$
 and  $d(k)g(k) \rightarrow 0$. Next we prove that 
 $d(k)\rightarrow 0$
 or  $g(k) \rightarrow 0$ are necessary conditions.
 
 \hfill $\bullet$

\vspace{0.5cm}

{\it Proposition 1}. We have the following two possibilities
\begin{equation}\label{P5}
or   \quad 
d(k) \rightarrow 0
\quad 
or \quad 
g(k) \rightarrow 0.
\end{equation}

{\it Proof.} We recall 
(\ref{B5A}). If
$c(k)\rightarrow \infty$
or
$h(k)\rightarrow \infty$,
then  
$\lim_{k\to\infty}
e^{2k(a-b)}c(k)h(k)$  it might not exist.
First, we prove here  that
both $c(k), h(k)$ have bounded limits limits to infinity.

Indeed, suppose
$c(k)\rightarrow \infty$.
That means, in fact, 
$E_a/\sigma(k) \rightarrow \infty$, because only this term of $c(k)$ is depending on $k$. 
We use the notation $I(k):=E_a(k)/\sigma$ and get 
\begin{equation}\label{B5C}
 e^{2k(a-b)}I^2(k) [\frac
 {\mu_L-\mu}{I(k)}-1]
 [\frac{\mu_R-\mu}{I(k)}-
\frac{E_b(k)}{E_a(k)}]
 -I^2(k)[\frac
 {\mu_L+\mu}{I(k)}-1]
 [\frac{\mu_R+\mu}{I(k)}-
\frac{E_b(k)}{E_a(k)}]=0.
 \end{equation}
 We simplify with $I^2(k)$. As
 $E_b(k)/E_a(k) 
 \rightarrow T_b/T_a$
 (which is a bounded quantity), 
 when $k\rightarrow\infty$
 we get 
$ 0 - T_a/T_b =0$.
Therefore the equation
 (\ref{B5A}) is not verified when 
 $c(k) \rightarrow \infty$.
 
 Suppose 
 $c(k), h(k) \rightarrow
 \infty$. Thus 
 $I(k):=E_a(k)\sigma(k)
  \rightarrow \infty, 
 \quad 
  J(k):=E_b(k)\sigma(k)
  \rightarrow \infty$ and from (\ref{B5A}) we get
\begin{equation}\label{B5E}
 e^{2k(a-b)}I(k)J(k) [\frac
 {\mu_L-\mu}{I(k)}-1]
 [\frac
 {\mu_R-\mu}{J(k)}-1]
 -I(k)J(k)
 [\frac
 {\mu_L+\mu}{I(k)}-1]
 [\frac
 {\mu_R+\mu}{J(k)}-1]=0.
\end{equation}
When $k \rightarrow \infty$ we get   
$0-1=0$.
Thus (\ref{B5A}) is not verified when 
 $c(k), h(k) \rightarrow \infty$.

The second step is following. 
Both 
$c(k), h(k)$ have bounded limits when $k\rightarrow \infty$, therefore 
$e^{2k(a-b)}c(k)h(k)
\rightarrow 0.$
From (\ref{B5A}) 
it follows 
$ d(k)g(k) \rightarrow 0$ and we obtain 
the result  (\ref{P5}). 

\hfill $\bullet$

\vspace{0.25cm}

{\it Remark 3.}
The above proposition gives us  two possible 
eigenvalues,  
{\it only}  for  large $k$:
\begin{equation}\label{B5F}
k \rightarrow \infty 
\Rightarrow
\sigma_1(k) \approx
E_a(k)/(\mu_L+\mu), \quad 
\sigma_2(k) \approx
E_b(k)/(\mu_R+\mu).
\end{equation}
This was also obtained in \cite{PA3}
by direct calculations
from (\ref{B5A}), with   
the natural hypothesis
$$ e^{-2k} k^3 \approx 0, \,\,
e^{-2k} << k^3;   \quad  
e^{-2k} k^6 \approx  0, 
\,\,
e^{-2k} k^6 << k^6 \quad 
\mbox{ for large enough}
\quad k.              $$

\hfill $\bullet$

 \vspace{0.25cm}
 
 {\it Proposition 2}.  
 There exists unbounded  eigenfunctions $f(x,k)$ 
  with $A(k), B(k)$ given by  (\ref{B3}).

 \vspace{0.25cm}
 
{\it Proof.} Consider  
 $d(k) \rightarrow 0$
 and $f(a,k)$ bounded. Then from 
$(\ref{B3})_1$ 
 we obtain
\begin{equation}\label{B11} (-2\mu)
 A(k)e^{ka} \rightarrow 0.
\end{equation}
Let $Y(k)=-g(k)/h(k)$. From (\ref{P1}),
(\ref{P2}),
(\ref{P3}) it follows that 
$Y(k)$ is bounded 
 for large $k$.
We use  $(\ref{B3})_2$
and  obtain
\begin{equation}\label{B12}
 B(k)e^{-ka} = Y(k) A(k)
 e^{k(2b-a)}.
 \end{equation}
 Here is the main point of our paper. We use
 (\ref{B12})
 to get the limit of 
 $B(k)e^{-ka}$.
There exist a lot of possibilities, related with  $A(k)$ and  the magnitude of the interval $(a,b)$. For example:

\vspace{0.2cm}

i) if $A(k)=D=constant$
and $a=-10, b=-1$, then
$2b-a=8>0$ and
$|B(k)e^{-ka}| \rightarrow 
\infty$. 
This contradicts
the hypothesis
$f(a,k)$ bounded;

\vspace{0.2cm}

ii) if $A=D$ and 
$a=-2, b=-1,5$ then
$2b-a= -1<0$ and both
terms of $f(a,k)$ are convergent to zero for  increasing $k$. This is  a ``good'' case.

\vspace{0.2cm}

iii)  if $A(k)$ is a polynomial, 
$a=-10, b=-1$, we still have
$ A(k)e^{ka} \rightarrow 0$  and $|B(k)e^{-ka}| \rightarrow \infty.  $
This contradicts again 
the hypothesis
$f(a,k)$ bounded.
\vspace{0.2cm}

Moreover, the  formulas 
(\ref{B11})-(\ref{B12}) give us
\begin{equation}
\label{B13}
\lim_{k\to\infty}
 B(k)e^{-ka} = 
\lim_{k\to\infty} 
\{ Y(k) 
 [A(k)e^{ka}]
 e^{2k(b-a)} \}.
\end{equation}
 Thus
 $\lim_{k\to\infty}f(a,k)$ might even not exist.

\hfill $\bullet$

\vspace{0.5cm}

{\it Remark 4}. Consider
the case iii) with $b=0$  and formula (\ref{B13}) 
with $A(k)=e^{-ka}$. We get  a 
 very strange conclusion.
{\it 
 The eigenfunctions (so the perturbations amplitudes)
 could be  exponentially increasing, for fixed $k$, 
as functions of $(b-a)$. 
}
This is in total contradiction with the results of
\cite{CP}, \cite{DP} and 
\cite{GOR-HOM-1}. 
In these three papers,  a variable profile was considered
in the intermediate region - denoted by (I.R).
On the water-liquid interface, it was considered a
continuous viscosity. The main conclusion was: a  large enough  (I.R.)
(and a suitable viscosity  profile) can almost suppress the instability
Saffman-Taylor. This  effect was  confirmed by  experimental results. 
On the contrary, there are no experimental results for the models studied in
\cite{DA0} and 
\cite{DA1} - \cite{DA6}.

\hfill $\bullet$

\newpage

\begin{center}
{\bf 4. Conclusions}
\end{center}

In the seminal paper \cite{DA0} were given the following results.
From the relations  
(\ref{A6})-(\ref{A7})
was obtained the   system 
(\ref{B3}) with  $A,B$  
constant. 
The condition (\ref{B5A}) was used to get 
a quadratic equation for  
$\sigma$. 
Numerical values and some estimates of  
$\sigma$ were also given.
All these results were used later in 
\cite{DA1} - \cite{DA6}.
Some weak points  exist
in these cited papers,
which were explained in
\cite{PA00} -  \cite{PA3}.

 In this paper we consider
 the general case when 
 $A,B$ are depending on $k$.
 The eigenfunctions $f(x,k)$ are not bounded (in general) as functions of $k$, so the  amplitudes 
 of perturbations could be not small. Thus the solution of the stability system has no physical meaning, even if
 $A=A(k), B=B(k)$. 
 The equation
(\ref{B5A}) is still true. But we can say that  the corresponding eigenvalues also do not make physical sense.



\end{document}